\documentclass{amsart}
\usepackage{amssymb} 
\usepackage{color}
\usepackage[all]{xy}
\usepackage{enumerate}
\usepackage{ulem}

 \newtheorem{theorem}{Theorem}[section]

 \newtheorem{proposition}[theorem]{Proposition}
 \newtheorem{lemma}[theorem]{Lemma}
 \newtheorem{corollary}[theorem]{Corollary}
 \newtheorem{remark}[theorem]{Remark}
 \newtheorem{definition}[theorem]{Definition}
 \newtheorem{notation}[theorem]{Notation}

 \numberwithin{equation}{section}
 \def\subrel#1#2{\mathrel{\mathop{#2}\limits_{#1}}}

\def\cat{{\rm Cat}_{\infty}}

\def\wcat{\widehat{\rm Cat}_{\infty}}

\begin{document}

\title
{
A
perfect pairing for monoidal adjunctions}
\author{Takeshi Torii}
\address{Department of Mathematics, 
Okayama University,
Okayama 700--8530, Japan}
\email{torii@math.okayama-u.ac.jp}



\date{February 3, 2023 (version~4.0)}

\begin{abstract}
We give another proof
of the fact that there is a dual
equivalence between the $\infty$-category
of monoidal $\infty$-categories
with left adjoint oplax monoidal functors 
and that with right adjoint lax monoidal functors
by constructing a perfect pairing between them.
\end{abstract}

\maketitle

\section{Introduction}

In category theory
the right adjoint of an oplax monoidal
functor between monoidal categories
is lax monoidal.
It is reasonable to expect the similar statement
holds in higher category theory.

In \cite{Lurie1,Lurie2}
Lurie has developed higher category theory
using quasi-categories,
which are models for $(\infty,1)$-categories.
In particular, he has proved that
the right adjoint of a strong monoidal
functor is lax monoidal
between monoidal $\infty$-categories
in \cite{Lurie2}.
In \cite{Torii2}
we showed that the right adjoint
of an oplax monoidal functor is lax monoidal.
Haugseng-Hebestreit-Linskens-Nuiten~\cite{HHLN1}
proved that
the $\infty$-category ${\rm Mon}_{\mathcal{O}}^{\rm oplax,L}
    ({\rm Cat}_{\infty})$
of $\mathcal{O}$-monoidal $\infty$-categories 
with left adjoint oplax $\mathcal{O}$-monoidal functors 
is dual equivalent to the $\infty$-category
${\rm Mon}_{\mathcal{O}}^{\rm lax, R}
    ({\rm Cat}_{\infty})$
of $\mathcal{O}$-monoidal $\infty$-category
with right adjoint lax $\mathcal{O}$-monoidal functors
for each $\infty$-operad $\mathcal{O}^{\otimes}$.

\begin{theorem}[{\cite{HHLN1}}]
\label{thm:main-theorem}
There is an equivalence
\[ {\rm Mon}_{\mathcal{O}}^{\rm oplax,L}
    ({\rm Cat}_{\infty})^{\rm op}
    \stackrel{\simeq}{\longrightarrow}
    {\rm Mon}_{\mathcal{O}}^{\rm lax, R}
    ({\rm Cat}_{\infty})\]
of $\infty$-categories,
which is identity on objects,
and which assigns to a left adjoint 
oplax $\mathcal{O}$-monoidal functor
its right adjoint lax $\mathcal{O}$-monoidal functor.
\end{theorem}

The purpose of this note is to give
another proof of Theorem~\ref{thm:main-theorem}.
For this purpose 
we study functorialities of 
the monoidal Yoneda embeddings.
After that,
we will prove 
Theorem~\ref{thm:main-theorem}
by constructing a perfect pairing 
between 
${\rm Mon}_{\mathcal{O}}^{\rm oplax,L}
    ({\rm Cat}_{\infty})$ 
and 
${\rm Mon}_{\mathcal{O}}^{\rm lax, R}
    ({\rm Cat}_{\infty})$
(Theorem~\ref{thm:main-this-note}).
The author 
thinks
that the method of the proof 
might be interesting in its own right.

In fact,
Haugseng-Hebestreit-Linskens-Nuiten~\cite{HHLN1}
proved a more general statement
than Theorem~\ref{thm:main-theorem}.
They showed that 
there is a bidual equivalence 
between the $(\infty,2)$-category
of $\mathcal{O}$-monoidal $\infty$-categories
with left adjoint oplax $\mathcal{O}$-monoidal functors
and that of 
$\mathcal{O}$-monoidal $\infty$-categories
with right adjoint lax $\mathcal{O}$-monoidal functors.
We will not give a proof of the $(\infty,2)$-categorical
statement since it needs a more sophisticated method
to organize equivalences between mapping
$(\infty,1)$-categories
of the two $(\infty,2)$-categories.   
Instead,
we will compare the equivalence in this note with
the restriction of the equivalence 
of Haugseng-Hebestreit-Linskens-Nuiten
in \cite{Torii5}.

Now, we will describe an outline of our proof. 
A small $\infty$-category $\mathcal{C}$
can be embedded into the $\infty$-category ${\rm P}(\mathcal{C})$
of presheaves on $\mathcal{C}$ by the Yoneda lemma.
If $\mathcal{C}$ is an $\mathcal{O}$-monoidal
$\infty$-category,
then we can construct an $\mathcal{O}$-monoidal
$\infty$-category $\mathbb{P}_{\mathcal{O}}(\mathcal{C})$
such that $\mathbb{P}_{\mathcal{O}}(\mathcal{C})_X
\simeq {\rm P}(\mathcal{C}_X)$ for each $X\in\mathcal{O}$
by Day convolution product~\cite{Lurie2}.
We will construct a pairing between 
the $\infty$-category
${\rm Mon}_{\mathcal{O}}^{\rm oplax}({\rm Cat}_{\infty})$
of $\mathcal{O}$-monoidal $\infty$-categories
with oplax $\mathcal{O}$-monoidal functors
and the $\infty$-category
${\rm Mon}_{\mathcal{O}}^{\rm lax}({\rm Cat}_{\infty})$
of $\mathcal{O}$-monoidal $\infty$-categories
with lax $\mathcal{O}$-monoidal functors.
The pairing corresponds to
a functor
${\rm Mon}_{\mathcal{O}}^{\rm oplax}({\rm Cat}_{\infty})
\times
{\rm Mon}_{\mathcal{O}}^{\rm lax}({\rm Cat}_{\infty})
\to\widehat{\mathcal{S}}$,
which assigns to
a pair $(\mathcal{C},\mathcal{D})$
the mapping space
${\rm Map}_{{\rm Mon}_{\mathcal{O}}^{\rm lax}({\rm Pr}^{\rm L})}
(\mathbb{P}_{\mathcal{O}}(\mathcal{C}),
\mathbb{P}_{\mathcal{O}}(\mathcal{D}))$.
By restricting the pairing
to the full subcategory spanned by those vertices
corresponding to
right adjoint lax $\mathcal{O}$-monoidal
functors $\mathcal{C}\to \mathcal{D}$,
we obtain a perfect pairing between
${\rm Mon}_{\mathcal{O}}^{\rm oplax, L}({\rm Cat}_{\infty})$
and
${\rm Mon}_{\mathcal{O}}^{\rm lax, R}({\rm Cat}_{\infty})$,
which gives the desired dual equivalence.

The organization of this note is as follows:
In \S\ref{section:monoidal-structure-on-presheaves}
we study monoidal structures
on $\infty$-categories of presheaves
and Yoneda embeddings.
In \S\ref{section:presheaf-functor}
we consider monoidal functorialities
for the construction of $\infty$-categories of presheaves.
In \S\ref{section:lax-oplax-duality}
we prove Theorem~\ref{thm:main-theorem}
by constructing the perfect pairing. 


\begin{notation}\rm
We fix Grothendieck universes 
$\mathcal{U}\in \mathcal{V}\in\mathcal{W}
\in \mathcal{X}$
throughout this note. 
We say that an element of $\mathcal{U}$ is small,
an element of $\mathcal{V}$ is large, 
an element of $\mathcal{W}$ is very large,
and an element of $\mathcal{X}$
is super large.

We denote by ${\rm Cat}_{\infty}$
the large $\infty$-category of small $\infty$-categories,
and by $\widehat{\rm Cat}_{\infty}$
the very large $\infty$-category of large $\infty$-categories.
We denote by $\mathrm{Pr}^{\rm L}$ 
the subcategory of $\widehat{\rm Cat}_{\infty}$
spanned by presentable $\infty$-categories
and left adjoint functors.
We write $\mathcal{S}$
for the large $\infty$-category of small spaces,
$\widehat{\mathcal{S}}$
for the very large $\infty$-category
of large spaces,
and $\widetilde{S}$
for the super large $\infty$-category
of very large spaces.

For a small $\infty$-category $\mathcal{C}$,
we denote by ${\rm P}(\mathcal{C})$
the $\infty$-category 
${\rm Fun}(\mathcal{C}^{\rm op},\mathcal{S})$
of presheaves on $\mathcal{C}$
with values in $\mathcal{S}$.
For an $\infty$-operad $\mathcal{O}^{\otimes}\to
{\rm Fin}_*$,
we write $\mathcal{O}$ for $\mathcal{O}^{\otimes}_{\langle 1\rangle}$.
For an $\mathcal{O}$-monoidal $\infty$-category
$\mathcal{C}^{\otimes}\to \mathcal{O}^{\otimes}$,
we denote by $\mathcal{C}$ 
the underlying $\infty$-category
$\mathcal{C}^{\otimes}\times_{\mathcal{O}^{\otimes}}\mathcal{O}$,
and we say that $\mathcal{C}$
is an $\mathcal{O}$-monoidal $\infty$-category
for simplicity.
\end{notation}

\noindent
{\bf Acknowledgements}.
The author would like to thank
Jonathan Beardsley
for letting him know 
references for the result by
Haugseng-Hebestreit-Linskens-Nuiten.
He would also like to thank the referee 
for useful comments and suggestions.
The author was partially supported 
by JSPS KAKENHI Grant Numbers JP17K05253.

\section{Monoidal structure on presheaves}
\label{section:monoidal-structure-on-presheaves}

Let $\mathcal{O}^{\otimes}$ be a small $\infty$-operad.
When $\mathcal{C}$ is a small $\mathcal{O}$-monoidal
$\infty$-category,
we can construct an $\mathcal{O}$-monoidal $\infty$-category
$\mathbb{P}_{\mathcal{O}}(\mathcal{C})^{\otimes}\to \mathcal{O}^{\otimes}$ 
by using Day convolution product.
In this section we study the $\mathcal{O}$-monoidal 
$\infty$-category 
$\mathbb{P}_{\mathcal{O}}(\mathcal{C})$.

First, we recall that 
there is an $\mathcal{O}$-monoidal
$\infty$-category
$\mathbb{P}_{\mathcal{O}}(\mathcal{C})^{\otimes}\to\mathcal{O}^{\otimes}$
by Day convolution product
if $\mathcal{C}$
is a small $\mathcal{O}$-monoidal $\infty$-category.
We denote by $\mathcal{S}^{\times}\to{\rm Fin}_*$
the symmetric monoidal $\infty$-category
for the $\infty$-category $\mathcal{S}$
of spaces with Cartesian
symmetric monoidal structure.
Taking pullback along the map 
$\mathcal{O}^{\otimes}\to{\rm Fin}_*$,
we obtain an $\mathcal{O}$-monoidal
$\infty$-category
$\mathcal{S}_{\mathcal{O}}^{\otimes}
\to\mathcal{O}^{\otimes}$,
where
$\mathcal{S}_{\mathcal{O}}^{\otimes}
=\mathcal{S}^{\times}\times_{{\rm Fin}_*}\mathcal{O}^{\otimes}$.
For an $\mathcal{O}$-monoidal $\infty$-category
$\mathcal{C}^{\otimes}\to\mathcal{O}^{\otimes}$, 
there is 
an opposite $\mathcal{O}$-monoidal
$\infty$-category $(\mathcal{C}^{\vee})^{\otimes}\to
\mathcal{O}^{\otimes}$ such that
$(\mathcal{C}^{\vee})^{\otimes}_X\simeq
\mathcal{C}_X^{\rm op}$
for each $X\in\mathcal{O}$.
By \cite[Construction~2.2.6.7]{Lurie2},
we can construct a fibration 
$\mathbb{P}_{\mathcal{O}}(\mathcal{C})^{\otimes}
\to\mathcal{O}^{\otimes}$
of $\infty$-operads,
where
$\mathbb{P}_{\mathcal{O}}(\mathcal{C})^{\otimes}
={\rm Fun}^{\mathcal{O}}
(\mathcal{C}^{\vee},\mathcal{S}_{\mathcal{O}})^{\otimes}$.
In fact, 
it is a coCartesian fibration
of $\infty$-operads
by \cite[Proposition~2.2.6.16]{Lurie2}.

\begin{lemma}[{\cite[Proposition~2.2.6.16]{Lurie2}}]
\label{lemma:presheaf-monoidal-Day-convolution}
If $\mathcal{C}$ is a small $\mathcal{O}$-monoidal
$\infty$-category,
then there is 
an $\mathcal{O}$-monoidal $\infty$-category
$\mathbb{P}_{\mathcal{O}}(\mathcal{C})^{\otimes}\to\mathcal{O}^{\otimes}$
in ${\rm Pr}^{\rm L}$
by Day convolution product
such that $\mathbb{P}_{\mathcal{O}}(\mathcal{C})^{\otimes}_X\simeq
{\rm P}(\mathcal{C}_X)$
for each $X\in\mathcal{O}$.
\end{lemma}

\begin{notation}\rm
For $\mathcal{O}$-monoidal $\infty$-categories
$\mathcal{C}$ and $\mathcal{D}$,
we denote by 
\[ {\rm Fun}_{\mathcal{O}}^{\rm lax}(\mathcal{C},\mathcal{D}) \]
the $\infty$-category ${\rm Alg}_{\mathcal{C}/\mathcal{O}}(\mathcal{D})$
of lax $\mathcal{O}$-monoidal functors
between $\mathcal{C}$ and $\mathcal{D}$.
\end{notation}


By the universal property of Day convolution product
(\cite[Definition~2.2.6.1]{Lurie2}),
we can characterize a lax $\mathcal{O}$-monoidal functor 
$\mathcal{C}\to \mathbb{P}_{\mathcal{O}}(\mathcal{D})$
by its associated functor
$\mathcal{C}\times_{\mathcal{O}}\mathcal{D}^{\vee}\to
\mathcal{S}_{\mathcal{O}}
$.

\begin{lemma}[{\cite[Definition~2.2.6.1]{Lurie2}}]
\label{lemma:adjoint-characterization-lax-functor}
Let $\mathcal{C}$ be a small
$\mathcal{O}$-monoidal $\infty$-category
and let $\mathcal{D}$ be 
an $\mathcal{O}$-monoidal $\infty$-category.
There is a lax $\mathcal{O}$-monoidal functor
$\mathbb{P}_{\mathcal{O}}(\mathcal{C})
\times_{\mathcal{O}}\mathcal{C}^{\vee}
\to 
\mathcal{S}_{\mathcal{O}}$
which induces an equivalence
\[ {\rm Fun}_{\mathcal{O}}^{\rm lax}
   (\mathcal{D},
   \mathbb{P}_{\mathcal{O}}(\mathcal{C}))
   \stackrel{\simeq}{\longrightarrow}
   {\rm Fun}_{\mathcal{O}}^{\rm lax}
   (\mathcal{D}\times_{\mathcal{O}}\mathcal{C}^{\vee},
   \mathcal{S}_{\mathcal{O}}) \]
of $\infty$-categories.
\end{lemma}

When $\mathcal{C}$ is a small $\infty$-category,
we have the Yoneda embedding
$j: \mathcal{C}\to{\rm P}(\mathcal{C})$
which is fully faithful,
and the $\infty$-category ${\rm P}(\mathcal{C})$
is freely generated by $j(\mathcal{C})$
under small colimits
in the sense of \cite[Theorem~5.1.5.6]{Lurie1}.
We will show that $j$ is promoted to 
a strong $\mathcal{O}$-monoidal
functor if $\mathcal{C}$ is a small $\mathcal{O}$-monoidal
$\infty$-category.

\begin{lemma}
\label{lemma:Yoneda-embedding-strong-monoidal}
Let $\mathcal{C}$ be a small $\mathcal{O}$-monoidal
$\infty$-category.
There is a strong $\mathcal{O}$-monoidal functor 
$J: \mathcal{C}\to \mathbb{P}_{\mathcal{O}}(\mathcal{C})$
such that
$J_X: \mathcal{C}_X\to \mathbb{P}_{\mathcal{O}}(\mathcal{C})_X\simeq
{\rm P}(\mathcal{C}_X)$
is 
equivalent to
the Yoneda embedding
for each $X\in\mathcal{O}$.
\end{lemma}

\proof
First, we shall show that
there is a lax $\mathcal{O}$-monoidal functor
$J: \mathcal{C}\to\mathbb{P}_{\mathcal{O}}(\mathcal{C})$
such that $J_X$ is 
equivalent to
the Yoneda embedding
for each $X\in\mathcal{O}$.
By Lemma~\ref{lemma:adjoint-characterization-lax-functor},
it suffices to construct  
a lax $\mathcal{O}$-monoidal functor
$\mathcal{C}^{\vee}\times_{\mathcal{O}}\mathcal{C}\to
\mathcal{S}_{\mathcal{O}}
$
whose
restriction
$\mathcal{C}_X^{\vee}\times\mathcal{C}_X
\to\mathcal{S}$
is equivalent to the mapping space functor
for each $X\in\mathcal{O}$.

By \cite[Example~5.2.2.23]{Lurie2},
we have a pairing of $\mathcal{O}$-monoidal
$\infty$-categories
\[ \lambda^{\otimes}:{\rm TwArr}(\mathcal{C})^{\otimes}
   \longrightarrow 
   \mathcal{C}^{\otimes}\times_{\mathcal{O}^{\otimes}}
   (\mathcal{C}^{\vee })^{\otimes}. \]
Note that this is a strong
$\mathcal{O}$-monoidal functor.
By taking opposite $\mathcal{O}$-monoidal $\infty$-categories,
we obtain a strong $\mathcal{O}$-monoidal
functor
\[ (\lambda^{\vee })^{\otimes}:
   ({\rm TwArr}(\mathcal{C})^{\vee })^{\otimes}
   \longrightarrow 
   (\mathcal{C}^{\vee})^{\otimes}\times_{\mathcal{O}^{\otimes}}
   \mathcal{C}^{\otimes}. \]
The restriction 
$(\lambda^{\vee})^{\otimes}_X:
({\rm TwArr}(\mathcal{C})^{\vee})^{\otimes}_X
   \to
   (\mathcal{C}^{\vee})^{\otimes}_X\times
   \mathcal{C}^{\otimes}_X$
is a left fibration for any $X\in\mathcal{O}^{\otimes}$.
By using coCartesian pushforward,
we can show that 
$(\lambda^{\vee})^{\otimes}$ is a left fibration.
Thus, we obtain a
$(\mathcal{C}^{\vee})^{\otimes}\times_{\mathcal{O}^{\otimes}}
\mathcal{C}^{\otimes}$-monoid object of $\mathcal{S}$.
This means 
there is a lax $\mathcal{O}$-monoidal
functor 
$\mathcal{C}^{\vee}\times_{\mathcal{O}}\mathcal{C}\to
\mathcal{S}_{\mathcal{O}}
$
by \cite[Proposition~2.4.2.5]{Lurie2}.
Furthermore,
we can see that the restriction 
over each $X\in\mathcal{O}$
is equivalent to the mapping space functor
by \cite[Proposition~5.2.1.11]{Lurie2}.

Next,
we shall show that 
$J$
is strong $\mathcal{O}$-monoidal.
Let $\phi: X\to Y$ be an active morphism
of $\mathcal{O}^{\otimes}$
with $Y\in\mathcal{O}$.
Let $C\in \mathcal{C}^{\otimes}_X$.
If $X\simeq X_1\oplus\cdots\oplus X_n$,
then $C\simeq C_1\oplus\cdots\oplus C_n$,
where $C_i\in \mathcal{C}_{X_i}$ for $1\le i\le n$.
There is a coCartesian morphism
$C\to C'$ of $\mathcal{C}^{\otimes}$ over $\phi$,
where $C'\simeq \bigotimes_{\phi}C_i$.
By \cite[Proof of Corollary~2.2.6.14]{Lurie2},
we have to show that 
\[ \alpha: \times\circ (\prod j(C_i))
           \longrightarrow
           j(C')\circ \bigotimes_{\phi} \]
exhibits $j(C')$ as a left Kan extension 
of $\times\circ \prod j(C_i)$
along $\bigotimes_{\phi}$:
\[ \xymatrix{
    \mathcal{C}_{X_1}^{\rm op}\times\cdots
    \mathcal{C}_{X_n}^{\rm op}
    \ar[r]^{\prod j(C_i)}
    \ar[d]_{\bigotimes_{\phi}}&
    \mathcal{S}\times\cdots\times
               \mathcal{S}\ar[r]^-{\times}&
   \mathcal{S}\\
  \mathcal{C}_{Y}^{\rm op}.\ar[rru]_{j(C')}& & \\ 
}\]
This follows from the facts that
a left Kan extension of a representable functor
is also representable by Yoneda's lemma
and that $\times\circ\prod j(C_i)$
is represented by the object
$(C_1,\ldots,C_n)
\in \mathcal{C}_{X_1}\times\cdots\times\mathcal{C}_{X_n}$.
\qed

\begin{remark}\label{remark:J-fully-faithful}
\rm
Note that
$J^{\otimes}: \mathcal{C}^{\otimes}\to 
\mathbb{P}_{\mathcal{O}}(\mathcal{C})^{\otimes}$ 
is fully faithful.
Hence,
it induces a fully faithful functor
\[ J_*: {\rm Fun}_{\mathcal{O}}^{\rm lax}
   (\mathcal{D},\mathcal{C})
   \longrightarrow
   {\rm Fun}_{\mathcal{O}}^{\rm lax}
   (\mathcal{D},\mathbb{P}_{\mathcal{O}}(\mathcal{C})) \]
for any small $\mathcal{O}$-monoidal
$\infty$-category $\mathcal{D}$.
We also call 
$J: \mathcal{C}\to 
\mathbb{P}_{\mathcal{O}}(\mathcal{C})$ 
the
Yoneda embedding.
\end{remark}

\begin{definition}\rm
Suppose that $\mathcal{C}$ and $\mathcal{D}$
are compatible with small colimits in the sense
of \cite[Variant~3.1.1.19]{Lurie2}.
We denote by
\[ {\rm Fun}^{\rm lax,cp}_{\mathcal{O}}
   (\mathcal{C},\mathcal{D}) \]
the full subcategory of 
${\rm Fun}^{\rm lax}_{\mathcal{O}}(\mathcal{C},\mathcal{D})$
spanned by those lax $\mathcal{O}$-monoidal
functors $f$ such that
$f_X: \mathcal{C}_X\to \mathcal{D}_X$ 
is colimit-preserving for each $X\in\mathcal{O}$.
\end{definition}

\begin{proposition}
\label{prop:presheaf-monoidal-universality}
Let $\mathcal{C}$ be a small $\mathcal{O}$-monoidal
$\infty$-category, and 
let $\mathcal{D}$ be an $\mathcal{O}$-monoidal
$\infty$-category which is compatible
with small colimits.
Then the Yoneda embedding
$J: \mathcal{C}\to\mathbb{P}_{\mathcal{O}}(\mathcal{C})$ induces an equivalence
\[ J^*: {\rm Fun}_{\mathcal{O}}^{\rm lax,cp}
   (\mathbb{P}_{\mathcal{O}}(\mathcal{C}),\mathcal{D})
   \stackrel{\simeq}{\longrightarrow}
   {\rm Fun}_{\mathcal{O}}^{\rm lax}(\mathcal{C},\mathcal{D}) \]
of $\infty$-categories.
\end{proposition}

In order to prove 
Proposition~\ref{prop:presheaf-monoidal-universality},
we need the following lemma.

\begin{lemma}\label{lemma:cofinal-lemma}
Let $\mathcal{C}$ be a small $\mathcal{O}$-monoidal $\infty$-category.
For any $X\in\mathcal{O}$ and $G\in {\rm P}(\mathcal{C}_X)$,
the inclusion map $(\mathcal{C}_X)_{/G}\hookrightarrow
(\mathcal{C}^{\otimes}_{\rm act})_{/G}$
is cofinal. 
\end{lemma}

\proof
For an object $\overline{D}=(D,J^{\otimes}(D)\to G)$
of $(\mathcal{C}^{\otimes}_{\rm act})_{/G}$,
we set $\mathcal{E}_{\overline{D}}=
(\mathcal{C}_X)_{/G}\times_{(\mathcal{C}^{\otimes}_{\rm act})_{/G}}
((\mathcal{C}^{\otimes}_{\rm act})_{/G})_{\overline{D}/}$.
By \cite[Theorem~4.1.3.1]{Lurie1},
it suffices to show that
$\mathcal{E}_{\overline{D}}$ 
is weakly contractible for any $\overline{D}$.
We shall show that $\mathcal{E}_{\overline{D}}$ 
has an initial object,
which implies that it is weakly contractible. 
Let $\phi$ be an active morphism in $\mathcal{O}^{\otimes}$
over which $J^{\otimes}(D)\to G$ lies.
We take a coCartesian morphism
$D\to \phi_!D$ in $\mathcal{C}^{\otimes}$
covering $\phi$.
Since
$J^{\otimes}: \mathcal{C}^{\otimes}\to 
\mathbb{P}_{\mathcal{O}}(\mathcal{C})^{\otimes}$
preserves coCartesian morphisms
by Lemma~\ref{lemma:Yoneda-embedding-strong-monoidal},
we can factor $J^{\otimes}(D)\to G$
as $J^{\otimes}(D)\to J^{\otimes}(\phi_!D)\to G$.
Then we can see that
$(D\to \phi_!D, J^{\otimes}(D)\to J^{\otimes}(\phi_!D)\to G)$
is an initial object of $\mathcal{E}_{\overline{D}}$.
\qed

\proof[Proof of Proposition~\ref{prop:presheaf-monoidal-universality}]
By \cite[Proposition~3.1.1.20 and Corollary~3.1.3.4]{Lurie2},
the functor
\[ J^*: {\rm Fun}_{\mathcal{O}}^{\rm lax}
   (\mathbb{P}_{\mathcal{O}}(\mathcal{C}),\mathcal{D})
   \to
   {\rm Fun}_{\mathcal{O}}^{\rm lax}(\mathcal{C},\mathcal{D}) \]
admits a left adjoint $J_!$.
Lemma~\ref{lemma:cofinal-lemma} implies that
the restriction
$J_!(f)_X: \mathbb{P}_{\mathcal{O}}(\mathcal{C})_X
             \to
             \mathcal{D}_X$
can be identified with 
a left Kan extension of $f_X$
along the Yoneda embedding
$j: \mathcal{C}_X\to {\rm P}(\mathcal{C}_X)$
for any $X\in\mathcal{O}$.
Hence $J_!(f)$ lands in the full subcategory
${\rm Fun}_{\mathcal{O}}^{\rm lax, cp}(\mathbb{P}_{\mathcal{O}}
(\mathcal{C}),\mathcal{D})$
and $J^*J_!(f)\simeq f$.
In particular,
we see that $J_!$ is fully faithful
and the essential image is 
${\rm Fun}_{\mathcal{O}}^{\rm lax, cp}
(\mathbb{P}_{\mathcal{O}}(\mathcal{C}),\mathcal{D})$.
This completes the proof.
\qed

\begin{definition}\rm
For $\mathcal{O}$-monoidal
$\infty$-categories $\mathcal{C}$ and $\mathcal{D}$,
we set 
\[ {\rm Fun}_{\mathcal{O}}^{\rm oplax}(\mathcal{C},\mathcal{D})
   ={\rm Fun}_{\mathcal{O}}^{\rm lax}
    (\mathcal{C}^{\vee},\mathcal{D}^{\vee})^{\rm op}, \]
and call it the $\infty$-category
of oplax $\mathcal{O}$-monoidal functors.
For 
$
f
\in {\rm Fun}_{\mathcal{O}}^{\rm oplax}(\mathcal{C},\mathcal{D})$,
we say that 
$f$ is left adjoint if
$f_X: \mathcal{C}_X\to \mathcal{D}_X$
is left adjoint for all $X\in\mathcal{O}$. 
We denote by
\[ {\rm Fun}^{\rm oplax,L}_{\mathcal{O}}(\mathcal{C},\mathcal{D}) \]
the full subcategory
of ${\rm Fun}^{\rm oplax}_{\mathcal{O}}(\mathcal{C},\mathcal{D})$
spanned by
those oplax $\mathcal{O}$-monoidal functors
which are left adjoint.

For $f
\in {\rm Fun}_{\mathcal{O}}^{\rm lax}(\mathcal{C},\mathcal{D})$,
we also say that $f$ is right adjoint if
$f_X: \mathcal{C}_X\to \mathcal{D}_X$
is right adjoint for all $X\in\mathcal{O}$. 
We denote by
\[ {\rm Fun}^{\rm lax,R}_{\mathcal{O}}(\mathcal{C},\mathcal{D}) \]
the full subcategory
of ${\rm Fun}^{\rm lax}_{\mathcal{O}}(\mathcal{C},\mathcal{D})$
spanned by
those lax $\mathcal{O}$-monoidal functors
which are right adjoint.
\end{definition}

Note that there is an equivalence
\[ {\rm Fun}_{\mathcal{O}}^{\rm oplax, L}
   (\mathcal{C},\mathcal{D})\simeq
   {\rm Fun}_{\mathcal{O}}^{\rm lax, R}
   (\mathcal{C}^{\vee},\mathcal{D}^{\vee})^{\rm op}\]
of $\infty$-categories
by definition.

\begin{definition}\label{def:P-O-upper-star}
\rm
Let $\mathcal{C}$ and $\mathcal{D}$ be
small $\mathcal{O}$-monoidal $\infty$-categories.
We define a functor
\[ \mathbb{P}_{\mathcal{O},(\mathcal{C},\mathcal{D})}^*:
   {\rm Fun}_{\mathcal{O}}^{\rm oplax}
   (\mathcal{C},\mathcal{D})^{\rm op}
   \longrightarrow
   {\rm Fun}_{\mathcal{O}}^{\rm lax,cp}
   (\mathbb{P}_{\mathcal{O}}(\mathcal{D}),
    \mathbb{P}_{\mathcal{O}}(\mathcal{C})) \]
by the following commutative diagram 
\begin{align}\label{align:definition-P-star}
    \begin{array}{ccccc}
    {\rm Fun}_{\mathcal{O}}^{\rm oplax}
    (\mathcal{C},\mathcal{D})^{\rm op}&
    \stackrel{\mathbb{P}_{\mathcal{O},(\mathcal{C},\mathcal{D})}^*}
             {\hbox to 10mm{\rightarrowfill}}
    &{\rm Fun}_{\mathcal{O}}^{\rm lax,cp}
    (\mathbb{P}_{\mathcal{O}}(\mathcal{D}),
     \mathbb{P}_{\mathcal{O}}(\mathcal{C}))
    &\stackrel{J^*}{\subrel{\simeq}{\hbox to 10mm{\rightarrowfill}}}&
    {\rm Fun}_{\mathcal{O}}^{\rm lax}
    (\mathcal{D},\mathbb{P}_{\mathcal{O}}(\mathcal{C}))\\
    \mbox{$\scriptstyle\simeq$}
   \bigg\downarrow
    \phantom{\mbox{$\scriptstyle\simeq$}}
    & & & & 
    \phantom{\mbox{$\scriptstyle\simeq\hspace{1mm}$}}
   \bigg\downarrow
    \mbox{$\scriptstyle\simeq\hspace{1mm}$}    \\
    {\rm Fun}_{\mathcal{O}}^{\rm lax}
    (\mathcal{C}^{\vee},\mathcal{D}^{\vee})
    &\stackrel{J_*}{\hbox to 10mm{\rightarrowfill}}&
    {\rm Fun}_{\mathcal{O}}^{\rm lax}
    (\mathcal{C}^{\vee},
     \mathbb{P}_{\mathcal{O}}(\mathcal{D}^{\vee}))&
    \stackrel{\simeq}{\hbox to 10mm{\rightarrowfill}}&
    {\rm Fun}_{\mathcal{O}}^{\rm lax}
    (\mathcal{C}^{\vee}\times_{\mathcal{O}}\mathcal{D},
    \mathcal{S}_{\mathcal{O}})\\
   \end{array}
\end{align}
in $\wcat$,
where the top right horizontal arrow is an equivalence
by Proposition~\ref{prop:presheaf-monoidal-universality}
and the right vertical and 
the bottom right horizontal arrows are equivalences
by Lemma~\ref{lemma:adjoint-characterization-lax-functor}.
Note that 
$\mathbb{P}_{\mathcal{O},(\mathcal{C},\mathcal{D})}^*$ is fully faithful
since $J_*$ is fully faithful
by Remark~\ref{remark:J-fully-faithful}.
For an oplax $\mathcal{O}$-monoidal
functor $f: \mathcal{C}\to \mathcal{D}$,
we simply write $f^*$ for 
$\mathbb{P}_{\mathcal{O},(\mathcal{C},\mathcal{D})}^*(f)$.
\end{definition}

\begin{proposition}
\label{prop:oplax-upper-star-lax-functor}
Let $f: \mathcal{C}\to \mathcal{D}$
be an oplax $\mathcal{O}$-monoidal functor
between small $\mathcal{O}$-monoidal
$\infty$-categories.
Then $(f^*)_X$ is equivalent to the functor
$(f_X)^*: {\rm P}(\mathcal{D}_X)\to {\rm P}(\mathcal{C}_X)$
for each $X\in\mathcal{O}$.
\end{proposition}

\proof
This follows by restricting
diagram~(\ref{align:definition-P-star}) over $X\in\mathcal{O}$.
\qed

\begin{corollary}
\label{cor:oplax-left-lax-right}
Let $f: \mathcal{C}\to \mathcal{D}$
be a 
left adjoint
oplax $\mathcal{O}$-monoidal functor
between small $\mathcal{O}$-monoidal $\infty$-categories.
Then
there exists a 
right adjoint
lax $\mathcal{O}$-monoidal functor
$f^R: \mathcal{D}\to \mathcal{C}$
such that $(f^R)_X\simeq (f_X)^R$,
where $(f_X)^R$ is a right adjoint
to $f_X$.
In this case
we have a commutative diagram
\[ \begin{array}{ccc}
    \mathcal{D}
    &\stackrel{f^R}{\longrightarrow}&
    \mathcal{C}\\
    \mbox{$\scriptstyle J$}
   \bigg\downarrow
    \phantom{\mbox{$\scriptstyle J$}}
    & &
    \phantom{\mbox{$\scriptstyle J$}}
   \bigg\downarrow
    \mbox{$\scriptstyle J$}\\
    \mathbb{P}_{\mathcal{O}}(\mathcal{D})
    &\stackrel{f^*}{\longrightarrow}&
    \mathbb{P}_{\mathcal{O}}(\mathcal{C})\\
   \end{array}\]
in ${\rm Mon}_{\mathcal{O}}^{\rm lax}(\wcat)$.
\end{corollary}

\proof
By Lemma~\ref{lemma:Yoneda-embedding-strong-monoidal}
and 
Definition~\ref{def:P-O-upper-star},
the composite $f^*\circ J:\mathcal{D}\to \mathbb{P}_{\mathcal{O}}(\mathcal{C})$
is lax $\mathcal{O}$-monoidal.
We can regard $\mathcal{C}$
as a full subcategory of $\mathbb{P}_{\mathcal{O}}(\mathcal{C})$
through $J$
by Remark~\ref{remark:J-fully-faithful}.
For each $X\in\mathcal{O}$,
we have a commutative diagram
\[ \begin{array}{ccc}
     \mathcal{D}_X & \stackrel{(f_X)^R}{\longrightarrow} &
     \mathcal{C}_X \\
     \mbox{$\scriptstyle j$}
    \bigg\downarrow
     \phantom{\mbox{$\scriptstyle j$}}
     & & 
     \phantom{\mbox{$\scriptstyle j$}}
    \bigg\downarrow
     \mbox{$\scriptstyle j$}\\
     \mathrm{P}(\mathcal{D}_X)&
     \stackrel{(f_X)^*}{\longrightarrow}&
     \mathrm{P}(\mathcal{C}_X),\\
  \end{array}\]
where the vertical arrows are the Yoneda embeddings.
By Lemma~\ref{lemma:Yoneda-embedding-strong-monoidal}
and 
Proposition~\ref{prop:oplax-upper-star-lax-functor},
$f^*\circ J$
factors through $J: \mathcal{C}\to
\mathbb{P}_{\mathcal{O}}(\mathcal{C})$,
and we obtain the desired lax $\mathcal{O}$-monoidal
functor $f^R: \mathcal{D}\to\mathcal{C}$.
\qed

\bigskip

\begin{remark}\rm
Let $f: \mathcal{C}\to\mathcal{D}$ be a right adjoint
lax $\mathcal{O}$-monoidal functor between
small $\mathcal{O}$-monoidal $\infty$-categories.
Then $f^{\vee}: \mathcal{C}^{\vee}\to\mathcal{D}^{\vee}$
is a left adjoint oplax $\mathcal{O}$-monoidal functor.
By Corollary~\ref{cor:oplax-left-lax-right},
we have a right adjoint lax $\mathcal{O}$-monoidal functor
$(f^{\vee})^R: \mathcal{D}^{\vee}\to\mathcal{C}^{\vee}$.
We set $f^L=((f^{\vee})^R)^{\vee}: \mathcal{D}\to \mathcal{C}$.
Then $f^L$ is a left adjoint oplax $\mathcal{O}$-monoidal functor
such that $(f^L)_X\simeq (f_X)^L$,
where $(f_X)^L: \mathcal{D}_X\to \mathcal{C}_X$
is a left adjoint to $f_X: \mathcal{C}_X\to\mathcal{D}_X$.
\end{remark}

\begin{definition}\label{def:P-O-shriek}
\rm
Let $\mathcal{C}$ and $\mathcal{D}$ be
small $\mathcal{O}$-monoidal $\infty$-categories.
We define a functor
\[ \mathbb{P}_{\mathcal{O}!,(\mathcal{C},\mathcal{D})}:
   {\rm Fun}_{\mathcal{O}}^{\rm lax}
   (\mathcal{C},\mathcal{D})
   \longrightarrow
   {\rm Fun}_{\mathcal{O}}^{\rm lax,cp}
   (\mathbb{P}_{\mathcal{O}}(\mathcal{C}),
    \mathbb{P}_{\mathcal{O}}(\mathcal{D})) \]
by the following commutative diagram 
\[ \xymatrix{
    {\rm Fun}_{\mathcal{O}}^{\rm lax}
    (\mathcal{C},\mathcal{D})
    \ar[rr]^{\mathbb{P}_{\mathcal{O}!,(\mathcal{C},\mathcal{D})}}
    \ar[dr]_{J_*}
    &&
    {\rm Fun}_{\mathcal{O}}^{\rm lax,cp}
    (\mathbb{P}_{\mathcal{O}}(\mathcal{C}),
     \mathbb{P}_{\mathbb{O}}(\mathcal{D}))
    \ar[dl]^{J^*}\\
    & {\rm Fun}_{\mathcal{O}}^{\rm lax}
       (\mathcal{C},
       \mathbb{P}_{\mathcal{O}}(\mathcal{D}))&\\
   }\]
in $\wcat$,
where the right vertical arrow is an equivalence
by Proposition~\ref{prop:presheaf-monoidal-universality}.
Note that 
$\mathbb{P}_{\mathcal{O}!,(\mathcal{C},\mathcal{D})}$ is fully faithful
since $J_*$ is fully faithful
by Remark~\ref{remark:J-fully-faithful}.
For a lax $\mathcal{O}$-monoidal
functor $f: \mathcal{C}\to \mathcal{D}$,
we simply write $f_!$ 
for $\mathbb{P}_{\mathcal{O}!,(\mathcal{C},\mathcal{D})}(f)$.
\end{definition}

\begin{remark}\label{remark:f-shriek-characterization}
\rm
Let $f: \mathcal{C}\to \mathcal{D}$
be a lax $\mathcal{O}$-monoidal functor
between small $\mathcal{O}$-monoidal $\infty$-categories.
By construction,
we have a commutative diagram
\[ \begin{array}{ccc}
     \mathcal{C}
     &\stackrel{f}{\longrightarrow}&
     \mathcal{D}\\
     \mbox{$\scriptstyle J$}
    \bigg\downarrow
     \phantom{\mbox{$\scriptstyle J$}}
     &&
     \phantom{\mbox{$\scriptstyle J$}}
    \bigg\downarrow
     \mbox{$\scriptstyle J$}\\
     \mathbb{P}_{\mathcal{O}}(\mathcal{C})
     &\stackrel{f_!}{\longrightarrow}&
     \mathbb{P}_{\mathcal{O}}(\mathcal{D})\\
   \end{array}\]
in ${\rm Mon}^{\rm lax}(\wcat)$.
\end{remark}

\begin{lemma}\label{lemma:star-R-shriek-R}
Let $\mathcal{C}$ and $\mathcal{D}$ be 
small $\mathcal{O}$-monoidal $\infty$-categories.
If $f: \mathcal{D}\to \mathcal{C}$
is a left adjoint oplax $\mathcal{O}$-monoidal functor,
then we have an equivalence
$f^* \simeq (f^R)_!$
in ${\rm Fun}_{\mathcal{O}}^{\rm lax,cp}
(\mathbb{P}_{\mathcal{O}}(\mathcal{C}),
\mathbb{P}_{\mathcal{O}}(\mathcal{D}))$.
Dually, if $g: \mathcal{C}\to\mathcal{D}$
is a right adjoint lax $\mathcal{O}$-monoidal functor,
then we have 
$(g^L)^*\simeq g_!$.
\end{lemma}

\proof
The Yoneda embedding 
induces an equivalence
of $\infty$-categories
${\rm Fun}_{\mathcal{O}}^{\rm lax,cp}
(\mathbb{P}_{\mathcal{O}}(\mathcal{C}),
\mathbb{P}_{\mathcal{O}}(\mathcal{D}))
\stackrel{\simeq}{\to}
{\rm Fun}_{\mathcal{O}}^{\rm lax}
(\mathcal{C},\mathbb{P}_{\mathcal{O}}(\mathcal{D}))$
by Proposition~\ref{prop:presheaf-monoidal-universality}.
If $f$ is a left adjoint oplax $\mathcal{O}$-monoidal functor,
then Corollary~\ref{cor:oplax-left-lax-right}
and Remark~\ref{remark:f-shriek-characterization}
imply that 
$f^*\circ J\simeq J\circ f^R\simeq (f^R)_!\circ J$
in
${\rm Fun}_{\mathcal{O}}^{\rm lax}
(\mathcal{C},\mathbb{P}_{\mathcal{O}}(\mathcal{D}))$.
Hence $f^*\simeq (f^R)_!$.
If $g$ is a right adjoint
lax $\mathcal{O}$-monoidal functor,
then we have $(g^L)^*\simeq
((g^L)^R)_!\simeq g_!$.
\qed

\begin{proposition}\label{prop:lax-oplax-duality-two-objects}
Let $\mathcal{C}$ and $\mathcal{D}$
be small $\mathcal{O}$-monoidal $\infty$-categories.
There is a natural equivalence
\[ 
   D_{(\mathcal{C},\mathcal{D})}
   : 
   {\rm Fun}_{\mathcal{O}}^{\rm oplax, L}
   (\mathcal{C},\mathcal{D})^{\rm op}
   \stackrel{\simeq}{\longrightarrow}
   {\rm Fun}_{\mathcal{O}}^{\rm lax,R}
   (\mathcal{D},\mathcal{C}) \]
of $\infty$-categories,
which makes the following diagram commute
\[ \xymatrix{
   {\rm Fun}_{\mathcal{O}}^{\rm oplax, L}
   (\mathcal{C},\mathcal{D})^{\rm op}
   \ar[rr]^{D_{(\mathcal{C},\mathcal{D})}}
   \ar[dr]_{\mathbb{P}_{\mathcal{O},(\mathcal{C},\mathcal{D})}^*}
   &&
   {\rm Fun}_{\mathcal{O}}^{\rm lax,R}
   (\mathcal{D},\mathcal{C})
   \ar[dl]^{\hspace{5mm}
            \mathbb{P}_{\mathcal{O}!,(\mathcal{C},\mathcal{D})}}\\
   &{\rm Fun}_{\mathcal{O}}^{\rm lax,cp}
    (\mathbb{P}_{\mathcal{O}}(\mathcal{D}),
     \mathbb{P}_{\mathcal{O}}(\mathcal{C})).&
}\]
in $\wcat$.
The functor $D_{(\mathcal{C},\mathcal{D})}$
associates to a left adjoint oplax 
$\mathcal{O}$-monoidal functor $f$
its right adjoint lax $\mathcal{O}$-monoidal functor
$f^R$.
\end{proposition}

\proof
Since $\mathbb{P}_{\mathcal{O},(\mathcal{C},\mathcal{D})}^*$ 
and
$\mathbb{P}_{\mathcal{O}!,(\mathcal{C},\mathcal{D})}$ 
are fully faithful,
it suffices to show that  
the essential images are the same
for the existence of 
the equivalence $D_{(\mathcal{C},\mathcal{D})}$.
This follows from Lemma~\ref{lemma:star-R-shriek-R}.
We also obtain the last part 
by Lemma~\ref{lemma:star-R-shriek-R}.
\qed

\section{Presheaf functors}
\label{section:presheaf-functor}

In this section we study
monoidal functorialities of the construction 
$\mathbb{P}_{\mathcal{O}}(\mathcal{C})$
from a small $\mathcal{O}$-monoidal
$\infty$-category $\mathcal{C}$.

\subsection{Contravariant
functor $\mathbb{P}_{\mathcal{O}}^{*}$}

There is a functor
\[ \mathrm{P}^*:\cat^{\rm op}\longrightarrow {\rm Pr}^{\rm L},\]
which associates to a small $\infty$-category $\mathcal{C}$
the $\infty$-category $\mathrm{P}(\mathcal{C})$
of presheaves on $\mathcal{C}$,
and to a functor $f:\mathcal{C}\to\mathcal{D}$
the functor $f^*: \mathrm{P}(\mathcal{D})
\to\mathrm{P}(\mathcal{C})$
that is obtained by composing with $f^{\rm op}$.
In this subsection we will construct
a functor 
\[ \mathbb{P}_{\mathcal{O}}^{*}:
   {\rm Mon}_{\mathcal{O}}^{
   \rm oplax
   }
   (\cat)^{\rm op}
   \longrightarrow
   {\rm Mon}_{\mathcal{O}}^{\rm lax}
   ({\rm Pr}^{\rm L}) 
\]    
which is a lifting of the functor $\mathrm{P}^*$.

First, we consider the composite of functors
\[ F_{\mathcal{O}}^*: {\rm Mon}_{\mathcal{O}}^{\rm oplax}(\cat)^{\rm op}
   \stackrel{\simeq}{\longrightarrow}
   {\rm Mon}_{\mathcal{O}}^{\rm lax}(\cat)^{\rm op}   
   \stackrel{}{\longrightarrow}
   {\rm Mon}_{\mathcal{O}}^{\rm lax}(\wcat), \]
where the first arrow is an equivalence
given by $\mathcal{C}\mapsto \mathcal{C}^{\vee}$
and the second arrow is given by
$\mathcal{C}\mapsto {\rm Fun}^{\mathcal{O}}
(\mathcal{C},\mathcal{S}_{\mathcal{O}})
=\mathbb{P}_{\mathcal{O}}(\mathcal{C}^{\vee})$.
We note that $F_{\mathcal{O}}^*$ assigns 
$\mathbb{P}_{\mathcal{O}}(\mathcal{C})$
to a small $\mathcal{O}$-monoidal
$\infty$-category $\mathcal{C}$.
For any oplax $\mathcal{O}$-monoidal functor
$f: \mathcal{C}\to \mathcal{D}$,
the restriction of 
$F_{\mathcal{O}}^*(f): \mathbb{P}_{\mathcal{O}}(\mathcal{D})
\to \mathbb{P}_{\mathcal{O}}(\mathcal{C})$
over $X\in\mathcal{O}$
is equivalent to 
$(f_X)^*: {\rm P}(\mathcal{D}_X)\to {\rm P}(\mathcal{C}_X)$.
Since $(f_X)^*$ admits a right adjoint $(f_X)_*$,
the functor $F_{\mathcal{O}}^*$
factors through ${\rm Mon}_{\mathcal{O}}^{\rm lax}({\rm Pr}^{\rm L})$.
Hence we obtain 
the desired functor
$\mathbb{P}_{\mathcal{O}}^*$

For small $\mathcal{O}$-monoidal
$\infty$-categories
$\mathcal{C}$ and $\mathcal{D}$,
the functor $\mathbb{P}_{\mathcal{O}}^*$
induces a functor
\[ {\rm Map}_{{\rm Mon}_{\mathcal{O}}^{\rm oplax}(\cat)}
   (\mathcal{C},\mathcal{D})
  \longrightarrow
   {\rm Map}_{{\rm Mon}_{\mathcal{O}}^{\rm lax}({\rm Pr}^{\rm L})}
   (\mathbb{P}_{\mathcal{O}}(\mathcal{D}),
    \mathbb{P}_{\mathcal{O}}(\mathcal{C})) \]
of mapping spaces.
We notice that it is equivalent to 
the functor obtained from 
$\mathbb{P}_{\mathcal{O},(\mathcal{C},\mathcal{D})}^*$
in Definition~\ref{def:P-O-upper-star}
by taking core.
Hence we obtain the following proposition.

\begin{proposition}
\label{prop:P-star-monoidal}
There is a functor
\[ \mathbb{P}_{\mathcal{O}}^{*}:
   {\rm Mon}_{\mathcal{O}}^{\rm oplax}(\cat)^{\rm op}
   \longrightarrow
   {\rm Mon}_{\mathcal{O}}^{\rm lax}({\rm Pr}^{\rm L}),\]  
which associates to a small $\mathcal{O}$-monoidal
$\infty$-category $\mathcal{C}$
the $\mathcal{O}$-monoidal $\infty$-category
$\mathbb{P}_{\mathcal{O}}(\mathcal{C})$,
and to an oplax $\mathcal{O}$-monoidal
functor $f: \mathcal{C}\to \mathcal{D}$
the lax $\mathcal{O}$-monoidal functor 
$f^*: \mathbb{P}_{\mathcal{O}}(\mathcal{D})\to \mathbb{P}_{\mathcal{O}}(\mathcal{C})$.
\end{proposition}

\subsection{Covariant
functor $\mathbb{P}_{\mathcal{O}!}$}

There is a functor
\[ \mathrm{P}_!: {\rm Cat}_{\infty}
   \longrightarrow {\rm Pr}^{\rm L},\]  
which associates to a small $\infty$-category
$\mathcal{C}$ the $\infty$-category $\mathrm{P}(\mathcal{C})$
of presheaves on $\mathcal{C}$,
and to a functor $f: \mathcal{C}\to \mathcal{D}$
the functor $f_!:\mathrm{P}(\mathcal{C})\to\mathrm{P}(\mathcal{D})$
that is obtained by left Kan extension. 
In this subsection we will construct
a functor 
\[ \mathbb{P}_{\mathcal{O}!}:
   {\rm Mon}_{\mathcal{O}}^{\rm lax}(\cat)
   \longrightarrow
   {\rm Mon}_{\mathcal{O}}^{\rm lax}
   ({\rm Pr}^{\rm L}) \]    
which is a lifting of $\mathrm{P}_!$.

\begin{remark}\rm
The functor ${\rm P}_!$ has two definitions:
One is as a free cocompletion and
the other is as a left adjoint to ${\rm P}^*$. 
They are recently shown to be equivalent 
in \cite{HHLN2}.
\end{remark}

\bigskip



First, we consider a functor
\[ F_{\mathcal{O}!}: {\rm Mon}_{\mathcal{O}}^{\rm lax}(\cat)^{\rm op}
   \longrightarrow
   {\rm Fun}({\rm Mon}_{\mathcal{O}}^{\rm lax}({\rm Pr}^{\rm L}),
   \widetilde{\mathcal{S}}
   ) \]
whose adjoint
${\rm Mon}_{\mathcal{O}}^{\rm lax}(\cat)^{\rm op}\times
{\rm Mon}_{\mathcal{O}}^{\rm lax}({\rm Pr}^{\rm L})\to
\widetilde{\mathcal{S}}
$
is given by
$(\mathcal{C},\mathcal{D})\mapsto
   {\rm Map}_{{\rm Mon}_{\mathcal{O}}^{\rm lax}(\wcat)}
   (\mathcal{C},\mathcal{D})$.
By Proposition~\ref{prop:presheaf-monoidal-universality},
the Yoneda embedding
induces an equivalence
$J^*:
   {\rm Map}_{{\rm Mon}_{\mathcal{O}}^{\rm lax}({\rm Pr}^{\rm L})}
   (\mathbb{P}_{\mathcal{O}}(\mathcal{C}),\mathcal{D})
   \stackrel{\simeq}{\to}
   {\rm Map}_{{\rm Mon}_{\mathcal{O}}^{\rm lax}(\wcat)}
   (\mathcal{C},\mathcal{D})$.
This implies that the functor 
$F_{\mathcal{O}!}$ factors through
the Yoneda embedding
${\rm Mon}_{\mathcal{O}}^{\rm lax}({\rm Pr}^{\rm L})^{\rm op}
   \to
   {\rm Fun}({\rm Mon}_{\mathcal{O}}^{\rm lax}({\rm Pr}^{\rm L}),
\widetilde{\mathcal{S}}
)$.
Hence we obtain the desired functor
$\mathbb{P}_{\mathcal{O}!}$

For small $\mathcal{O}$-monoidal
$\infty$-categories
$\mathcal{C}$ and $\mathcal{D}$,
the functor $\mathbb{P}_{\mathcal{O}!}$
induces a functor
\[ {\rm Map}_{{\rm Mon}_{\mathcal{O}}^{\rm lax}(\cat)}
   (\mathcal{C},\mathcal{D})
  \longrightarrow
   {\rm Map}_{{\rm Mon}_{\mathcal{O}}^{\rm lax}({\rm Pr}^{\rm L})}
   (\mathbb{P}_{\mathcal{O}}(\mathcal{C}),
    \mathbb{P}_{\mathcal{O}}(\mathcal{D})) \]
of mapping spaces.
We notice that it is equivalent to 
the functor obtained from
$\mathbb{P}_{\mathcal{O}!,(\mathcal{C},\mathcal{D})}$
in Definition~\ref{def:P-O-shriek}
by taking core.
Hence we obtain the following proposition.

\begin{proposition}\label{prop:monoidal-shriek-map}
There is a functor
\[ \mathbb{P}_{\mathcal{O}!}:
   {\rm Mon}_{\mathcal{O}}^{\rm lax}(\cat)
   \longrightarrow
   {\rm Mon}_{\mathcal{O}}^{\rm lax}({\rm Pr}^{\rm L}),\]
which associates to a small $\mathcal{O}$-monoidal
$\infty$-category $\mathcal{C}$
the $\mathcal{O}$-monoidal $\infty$-category
$\mathbb{P}_{\mathcal{O}}(\mathcal{C})$, 
and to 
a
lax $\mathcal{O}$-monoidal
functor $f: \mathcal{C}\to \mathcal{D}$
the lax $\mathcal{O}$-monoidal functor 
$f_!: \mathbb{P}_{\mathcal{O}}(
\mathcal{C}
)\to 
\mathbb{P}_{\mathcal{O}}(
\mathcal{D}
)$.
\end{proposition}

\section{Perfect pairing
for monoidal adjunctions}
\label{section:lax-oplax-duality}

In this section we will 
prove Theorem~\ref{thm:main-theorem}
by constructing a perfect pairing
between the $\infty$-category
of $\mathcal{O}$-monoidal $\infty$-categories
with left adjoint oplax monoidal
functors and that
with right adjoint lax monoidal functors
(Theorem~\ref{thm:main-this-note}).


We denote by
\[ \widehat{\lambda}:
   \widehat{\mathcal{M}}_{\mathcal{O}}^{\rm lax}
   \longrightarrow
   {\rm Mon}_{\mathcal{O}}^{\rm lax}({\rm Pr}^{\rm L})
   \times 
   {\rm Mon}_{\mathcal{O}}^{\rm lax}({\rm Pr}^{\rm L})^{\rm op} \]
the perfect pairing associated
to the mapping space functor
${\rm Map}_{{\rm Mon}_{\mathcal{O}}^{\rm lax}({\rm Pr}^{\rm L})}(-,-):
   {\rm Mon}_{\mathcal{O}}^{\rm lax}({\rm Pr}^{\rm L})^{\rm op}
   \times 
   {\rm Mon}_{\mathcal{O}}^{\rm lax}({\rm Pr}^{\rm L})
   \to\widehat{\mathcal{S}}$.
We define 
\[ \lambda:
   \mathcal{M}_{\mathcal{O}}^{\rm lax}\longrightarrow
   {\rm Mon}_{\mathcal{O}}^{\rm oplax,L}(\cat)^{\rm op}\times
   {\rm Mon}_{\mathcal{O}}^{\rm lax,R}(\cat)^{\rm op}
\]
to be a right fibration
obtained from $\widehat{\lambda}$
by pullback along the functor
\[ \mathbb{P}_{\mathcal{O}}^*\times
   \mathbb{P}_{\mathcal{O}!}^{\,\rm op}:
    {\rm Mon}_{\mathcal{O}}^{\rm oplax,L}(\cat)^{\rm op}\times
   {\rm Mon}_{\mathcal{O}}^{\rm lax,R}(\cat)^{\rm op}
   \longrightarrow 
   {\rm Mon}_{\mathcal{O}}^{\rm lax}({\rm Pr}^{\rm L})
   \times
   {\rm Mon}_{\mathcal{O}}^{\rm lax}({\rm Pr}^{\rm L})^{\rm op}. \]
An object of $\mathcal{M}^{\rm lax}_{\mathcal{O}}$
corresponds to a triple $(\mathcal{C},\mathcal{D},f)$,
where $\mathcal{C}$ and $\mathcal{D}$ are
small $\mathcal{O}$-monoidal $\infty$-categories
and $f: \mathbb{P}_{\mathcal{O}}(\mathcal{C})\to
\mathbb{P}_{\mathcal{O}}(\mathcal{D})$
is a morphism in 
${\rm Mon}_{\mathcal{O}}^{\rm lax}({\rm Pr}^{\rm L})$.

For small $\mathcal{O}$-monoidal
$\infty$-categories $\mathcal{C}$ and $\mathcal{D}$,
the functor $\mathbb{P}_{\mathcal{O}!}$
induces a fully faithful functor 
${\rm Map}_{{\rm Mon}_{\mathcal{O}}^{\rm lax}(\cat)}
   (\mathcal{C},\mathcal{D})
  \to
   {\rm Map}_{{\rm Mon}_{\mathcal{O}}^{\rm lax}({\rm Pr}^{\rm L})}
   (\mathbb{P}_{\mathcal{O}}(\mathcal{C}),
\mathbb{P}_{\mathcal{O}}(\mathcal{D}))$
of mapping spaces.
We define
\[ \mathcal{M}_{\mathcal{O}}^{\rm lax, R} \]
to be the full subcategory of 
$\mathcal{M}_{\mathcal{O}}^{\rm lax}$
spanned by those 
objects
$v$ 
corresponding to triples $(\mathcal{C},\mathcal{D},f)$
where $f\simeq g_!$
for some right adjoint lax $\mathcal{O}$-monoidal
functor $g: \mathcal{C}\to\mathcal{D}$.
We let 
\[ \lambda^{\rm R}: 
   \mathcal{M}_{\mathcal{O}}^{\rm lax, R}\longrightarrow
   {\rm Mon}_{\mathcal{O}}^{\rm oplax,L}(\cat)^{\rm op}\times
   {\rm Mon}_{\mathcal{O}}^{\rm lax,R}(\cat)^{\rm op}
\]
be the restriction
of $
\lambda
$ to 
$\mathcal{M}_{\mathcal{O}}^{\rm lax, R}$.

First,
we show that $\lambda^{\rm R}$ is a paring
of $\infty$-categories.

\begin{lemma}\label{lemma:lambda-R-pairing}
The functor
$\lambda^{\rm R}:
   \mathcal{M}_{\mathcal{O}}^{\rm lax, R}
   \to
   {\rm Mon}_{\mathcal{O}}^{\rm oplax,L}(\cat)^{\rm op}
   \times
   {\rm Mon}_{\mathcal{O}}^{\rm lax,R}(\cat)^{\rm op}$
is a pairing of $\infty$-categories. 
\end{lemma}

\proof
We shall prove that 
$\lambda^{\rm R}$ is a right fibration.
For this purpose,
since
$
\lambda
$ is a right fibration,
it suffices to show the following:
For a morphism $m: v\to v'$ of 
$
\mathcal{M}_{\mathcal{O}}^{\rm lax}
$,
if $v'\in \mathcal{M}_{\mathcal{O}}^{\rm lax, R}$, 
then $v\in\mathcal{M}_{\mathcal{O}}^{\rm lax, R}$.

We denote by 
$k: \mathbb{P}_{\mathcal{O}}(\mathcal{C})\to 
\mathbb{P}_{\mathcal{O}}(\mathcal{D})$ 
and
$k': \mathbb{P}_{\mathcal{O}}(\mathcal{C}')\to 
\mathbb{P}_{\mathcal{O}}(\mathcal{D}')$
morphisms in ${\rm Mon}_{\mathcal{O}}^{\rm lax}({\rm Pr}^{\rm L})$ 
corresponding to $v$ and $v'$, respectively.
Since $v'\in \mathcal{M}^{\rm lax, R}_{
\mathcal{O}
}$,
we can write $k'\simeq f_!$,
where $f: \mathcal{C}'\to\mathcal{D}'$
is a right adjoint lax $\mathcal{O}$-monoidal functor.
Let 
$
\lambda
(m)\simeq (g,h)$.
Then $k\simeq h_!\circ k'\circ g^*$.
Since $g$ is a left adjoint oplax 
$\mathcal{O}$-monoidal functor,
$g^*\simeq (g^R)_!$
by Lemma~\ref{lemma:star-R-shriek-R}.
Thus, 
$k\simeq (h\circ f\circ g^R)_!$.
Since $h\circ f\circ g^R$ 
is a right adjoint lax $\mathcal{O}$-monoidal
functor,
we see that 
$v\in\mathcal{M}^{\rm lax, R}_{
\mathcal{O}
}$.
\qed

\bigskip

Next, we shall show that
the pairing $\lambda^R$ is perfect.
For this purpose,
we need the following lemma.

\begin{lemma}\label{lemma:representing-equivalence-conditions}
Let $v$ be an 
object of $\mathcal{M}_{\mathcal{O}}^{\rm lax,R}$
which corresponds to
a triple $(\mathcal{C},\mathcal{D},f)$.
Then the following conditions are equivalent\,\mbox{\rm :}
\begin{enumerate}

\item[{\rm (1)}]
The object $v$ is left universal.

\item[{\rm (2)}]
The object $v$ is right universal.

\item[{\rm (3)}]
$f\simeq g_!$,
where $g: \mathcal{C}\to\mathcal{D}$
is an equivalence 
in ${\rm Mon}_{\mathcal{O}}^{\rm lax,R}(\cat)$.

\end{enumerate}
\end{lemma}

\proof
First, 
we shall prove the equivalence between (1) and (3).
Let 
$\lambda^{\rm R}_{\mathcal{C}}$
be the right fibration 
obtained from $\lambda^{\rm R}$
by restriction to
$\{\mathcal{C}\}
\times{{\rm Mon}_{\mathcal{O}}^{\rm lax,R}(\cat)^{\rm op}}$.
It is classified by
the functor 
${\rm Mon}_{\mathcal{O}}^{\rm lax, R}(\cat)\to \mathcal{S}$
given by
$\mathcal{D}'\mapsto
{\rm Map}_{{\rm Mon}_{\mathcal{O}}^{\rm lax,R}(\cat)}
(\mathcal{C},\mathcal{D}')$.
Hence 
$\lambda^{\rm R}_{\mathcal{C}}$
is equivalent to
the map
${\rm Mon}_{\mathcal{O}}^{\rm lax,R}(\cat)^{\rm op}_{/\mathcal{C}}
\to {\rm Mon}_{\mathcal{O}}^{\rm lax,R}(\cat)^{\rm op}$
as right fibrations.
Thus,
we see that 
$v$ is left universal if and only if 
$f\simeq g_!$ for an equivalence $g$.

Next,
we shall prove
the equivalence between (2) and (3).
We consider the right fibration
$\lambda^{\rm R}_{\mathcal{D}}:
\mathcal{N}\to 
{\rm Mon}_{\mathcal{O}}^{\rm oplax, L}(\cat)^{\rm op}$
obtained from $\lambda^{\rm R}$
by restriction to
${\rm Mon}_{\mathcal{O}}^{\rm oplax, L}(\cat)^{\rm op}
\times\{\mathcal{D}\}$,
where 
$\mathcal{N}=\mathcal{M}^{\rm lax, R}_{
\mathcal{O}
}
\times_{{\rm Mon}_{\mathcal{O}}^{\rm lax,R}(\cat)^{\rm op}}\{\mathcal{D}\}$.
%
We have
a commutative diagram
of right fibrations
\[ \begin{array}{ccc}
    \mathcal{N}_{/
    v
     } & \longrightarrow &
    \mathcal{N}\\[1mm]
   \bigg\downarrow    & & 
   \bigg\downarrow\\[4mm] 
    {\rm Mon}_{\mathcal{O}}^{\rm oplax, L}(\cat)^{\rm op}
    {}_{/
    \mathcal{C}
    }&
    \longrightarrow &
    {\rm Mon}_{\mathcal{O}}^{\rm oplax, L}(\cat)^{\rm op}.\\
   \end{array}\]
We notice that the left vertical arrow is 
an equivalence
since it is induced on overcategories by a right fibration.
We would like to show that
the map $\mathcal{N}_{/v}\to\mathcal{N}$
is an equivalence
if and only if $f$ is an equivalence.

For any $\mathcal{C}'\in 
{\rm Mon}_{\mathcal{O}}^{\rm oplax, L}(\cat)$,
we have a pullback diagram
\[ \begin{array}{ccc} 
   \mathcal{N}_{/v}\times_{{\rm Mon}_{\mathcal{O}}^{\rm oplax,L}(\cat)^{\rm op}}
   \{\mathcal{C}'\}
   &{\hbox to 15mm{\rightarrowfill}}&
   {\rm Map}_{{\rm Mon}_{\mathcal{O}}^{\rm lax,R}(\cat)}
   (\mathcal{C}',\mathcal{D})\\[2mm]
  \bigg\downarrow 
   & & 
   \phantom{\mbox{$\scriptstyle \mathbb{P}_{\mathcal{O}!}$}}
  \bigg\downarrow
   \mbox{$\scriptstyle \mathbb{P}_{\mathcal{O}!}$}\\
   {\rm Map}_{{\rm Mon}_{\mathcal{O}}^{\rm oplax,L}(\cat)^{\rm op}}
   (\mathcal{C}',\mathcal{C}) 
   & \stackrel{f\circ\mathbb{P}_{\mathcal{O}}^*}
   {\hbox to 15mm{\rightarrowfill}}&
   {\rm Map}_{{\rm Mon}_{\mathcal{O}}^{\rm lax}({\rm Pr}^{\rm L})}
   (\mathbb{P}_{\mathcal{O}}(\mathcal{C}'),
    \mathbb{P}_{\mathcal{O}}(\mathcal{D})) \\
   \end{array}   \]
in spaces,
where the left vertical arrow is an equivalence.
Suppose that $f\simeq g_!$
with $g:\mathcal{C}\to \mathcal{D}$
a right adjoint lax $\mathcal{O}$-monoidal functor.
By Lemma~\ref{lemma:star-R-shriek-R}
and
Proposition~\ref{prop:lax-oplax-duality-two-objects},
we see that the map
\[ \mathcal{N}_{/v}\times_{{\rm Mon}_{\mathcal{O}}^{\rm oplax,L}(\cat)^{\rm op}}
   \{\mathcal{C}'\}
   \longrightarrow
   \mathcal{N}\times_{{\rm Mon}_{\mathcal{O}}^{\rm oplax,L}(\cat)^{\rm op}}
   \{\mathcal{C}'\} \]
is equivalent to the composite
\[ {\rm Map}_{{\rm Mon}_{\mathcal{O}}^{\rm oplax, L}(\cat)}
   (\mathcal{C},
   \mathcal{C}'
   )
   \stackrel{\simeq}{\to}
   {\rm Map}_{{\rm Mon}_{\mathcal{O}}^{\rm lax, R}(\cat)}
   (
    \mathcal{C}',
   \mathcal{C})
   \stackrel{g\circ(-)}
   {\hbox to 12mm{\rightarrowfill}}
   {\rm Map}_{{\rm Mon}_{\mathcal{O}}^{\rm lax, R}(\cat)}
   (
    \mathcal{C}',
    \mathcal{D}), \]
where the first arrow is an equivalence induced by
Proposition~\ref{prop:lax-oplax-duality-two-objects}
and the second arrow is 
the composition with $g$. 
Thus, each fiber of the map
$\mathcal{N}_{/
v
}\to\mathcal{N}$
is contractible 
if and only if $g$ is an equivalence.
Therefore, 
(2) and (3) are equivalent.
\qed

\begin{proposition}\label{prop:M^R^perfect-pairing}
The paring $\lambda^{\rm R}$
is perfect.
\end{proposition}

\proof
The proposition follows
from \cite[Corollary~5.2.1.22]{Lurie2}
and Lemma~\ref{lemma:representing-equivalence-conditions}.
\qed

\bigskip

We obtain 
the main theorem of this note.

\begin{theorem}\label{thm:main-this-note}
There is an equivalence
\[ T: {\rm Mon}_{\mathcal{O}}^{\rm lax,R}(\cat)
      \stackrel{\simeq}{\longrightarrow}
      {\rm Mon}_{\mathcal{O}}^{\rm oplax,L}(\cat)^{\rm op} \]
of $\infty$-categories,
which is identity on objects and assigns
to right adjoint lax $\mathcal{O}$-monoidal functors
their left adjoint oplax $\mathcal{O}$-monoidal functors.
The equivalence $T$ fits into
the following commutative diagram
\[ \xymatrix{
    {\rm Mon}_{\mathcal{O}}^{\rm lax,R}(\cat)
    \ar[rr]^T\ar[dr]_{\mathbb{P}_{\mathcal{O}!}}
   &&
      {\rm Mon}_{\mathcal{O}}^{\rm oplax,L}(\cat)^{\rm op} 
    \ar[dl]^{\mathbb{P}_{\mathcal{O}}^*}\\
    &{\rm Mon}_{\mathcal{O}}^{\rm lax}({\rm Pr}^{\rm L}).&
   }\]
\end{theorem}

\proof
The first part follows from
Proposition~\ref{prop:M^R^perfect-pairing}.
We have a left (and right) representable morphism 
$\lambda^R\to \widehat{\lambda}$ of perfect pairings.
The second part follows from \cite[Proposition~5.2.1.17]{Lurie2}.
\qed

\end{document}